\documentclass[a4paper,11pt,reqno]{amsart}

\usepackage{amsmath}
\usepackage{amssymb}
\usepackage{amsfonts}
\usepackage[all,2cell]{xy} \UseAllTwocells \SilentMatrices
\usepackage{geometry,latexsym,amsmath,amstext,mathrsfs}
\usepackage{amssymb,amsthm,amscd,amsopn,verbatim,graphics}
\usepackage[english]{babel}

\theoremstyle{plain}
\newtheorem{theorem}{Theorem}
\newtheorem{lemma}[theorem]{Lemma}
\newtheorem{proposition}[theorem]{Proposition}

\theoremstyle{definition}
\newtheorem{definition}[theorem]{Definition}
\theoremstyle{remark}
\newtheorem{example}[theorem]{Example}
\newtheorem{remark}[theorem]{Remark}

\numberwithin{theorem}{section} \numberwithin{equation}{section}

\newcommand*\Hom{\mathrm{Hom}}

\def\P2{{P^{[2]}}}

\def\11{{\mbox{\boldmath $1$}}}

\def\eq{\begin{equation}}
\def\en{\end{equation}}

\def\eq#1\en{\begin{equation}#1\end{equation}}
\def\eqa#1\ena{\begin{eqnarray}#1\end{eqnarray}}

\begin{document}

\baselineskip=1.2\baselineskip

\title[From simplicial Lie algebras
to DGLA{\Tiny s} via 1-jets]{From simplicial Lie algebras and
hypercrossed complexes to differential graded Lie algebras via
1-jets}
\author{Branislav Jur\v
co}
\address{Mathematical Institute, Charles university, Sokolovsk\'a 83, Prague
186~75, Czech Republic} \email{jurco@karlin.mff.cuni.cz,
branislav.jurco@googlemail.com}
\begin{abstract}
Let $\mathfrak{g}$ be a simplicial Lie algebra with Moore complex
$N\mathfrak{g}$ of length $k$. Let $G$ be the simplicial Lie group
integrating $\mathfrak{g}$, such that each $G_n$ is simply
connected. We use the 1-jet of the classifying space $\overline WG$
to construct, starting from $\mathfrak{g}$, a Lie $k$-algebra $L$.
The so constructed Lie $k$-algebra $L$ is actually a differential
graded Lie algebra. The differential and the brackets are explicitly
described in terms (of a part) of the corresponding $k$-hypercrossed
complex structure of $N\mathfrak{g}$. The result can be seen as a
geometric interpretation of Quillen's (purely algebraic)
construction of the adjunction between simplicial Lie algebras and
dg-Lie algebras.
\end{abstract}

\maketitle

\section{Introduction}
In this paper we describe a geometric construction leading to
Quillen's relation between simplicial Lie algebras and differential
graded Lie algebras (DGLAs) \cite{Quillen}. We do that following the
ideas of \v Severa \cite{Sev}, which lead to a construction of
$L_\infty$-algebras (or, more generally, $L_\infty$-algebroids) as
1-jets (differentiation) of simplicial manifolds. Here, we apply \v
Severa's construction to the case when the simplicial manifold in
question is the classifying space $\overline{W}G$ of a simplicial
Lie group $G$, the simplicial Lie group $G$ having Moore complex of
length $k$. Main results are the Proposition \ref{da} and the
Theorem \ref{theMain}. In Proposition \ref{da} we describe
explicitly the dg manifold representing the 1-jet functor
$F_1^{\overline WG}$ and in Theorem \ref{theMain} we describe
explicitly  the corresponding $L_\infty$-algebra as a $k$-term
differential graded Lie algebra $L$ with the differential and
brackets given in terms the hypercrossed complex structure of
$N\mathfrak{g}$. The result is the same as the one described by the
$N$-functor in the Quillen's adjunction between simplicial Lie
algebras and dg-Lie algebras (see Proposition 4.4 of \cite{Quillen}.
The construction can equivalently be viewed as an assignment of a
$k$-term DGLA to a k-hypercrossed complex
$\underline{\mathfrak{g}}$. The paper is organized as follows.

In Section 2, we recall the relevant material concerning simplicial
Lie groups. In particular, we describe the Moore complex of a
simplicial Lie group and illustrate its Lie hypercrossed complex
structure in the low dimensional case of Lie crossed modules an Lie
2-crossed modules.

In Section 3, we recall the relevant facts regarding simplicial
principal bundles.

In Section 4, we summarize \v Severa's construction  and give the
relevant examples following \cite{Sev}. In particular, we describe
in detail the construction of a Lie algebra $\mathfrak{g}$ as a
1-jet of the classifying space $BG$ of the corresponding Lie group
$G$. Also, we describe in detail the construction of a Lie
$2$-algebra corresponding to a crossed module of Lie algebras
$\mathfrak{h}\to \mathfrak{d}$ as a 1-jet of the functor associating
to a surjective submersion $M\to N$ of (super)manifolds the set of
$(H\to D)$-descent data over $M\to N$.

The examples mentioned above are the starting point of this paper.
For this, we note that also the second example can be reinterpreted
as the 1-jet of a simplicial manifold. The relevant simplicial
manifold is the Duskin nerve of the strict Lie 2-group defined by
the Lie crossed module $H\to D$, which is isomorphic to $\overline W
G$, the classifying space of the simplicial Lie group associated to
the Lie crossed module $H\to D$. Therefore, it is natural to
generalize the above examples by applying \v Severa's construction
to the case of any simplicial Lie group $G$ and describe explicitly
the corresponding 1-jet of $\overline W G$. This is done in Section
5. The resulting dg manifold is described in Proposition \ref{da}
and the corresponding DGLA in Proposition \ref{theMain}. This DGLA
is the same as the one described by Quillen in Sec. 4 of
\cite{Quillen}.

In this paper we do not discuss, up to occasional
remarks,\footnote{Cf. Remarks \ref{gauge1} and \ref{gauge2}.}
applications to the higher gauge theory. These will be given in a
forthcoming paper.

All commutators are implicitly assumed to be graded. Although we do
not mention it explicitly, all constructions extend more or less
straightforwardly to the case when all in involved Lie groups and
Lie algebras are super. Hopefully, this is a wormless paper
\cite{KochanSev}.
\section{Simplicial groups and higher crossed modules}
Here we briefly sketch the relation between simplicial groups and
hypercrossed complexes of groups. The basic idea comes from
\cite{Cond} and is further developed and formalized in \cite{Car}.
We follow \cite{Cond, MutPor, Eh}.

Although the above mentioned references\footnote{and also our basic
reference regarding simplicial objects \cite{May} as well as other
useful references \cite{Curtis, Goerss}} work with simplicial sets,
the constructions and statements relevant relevant for our purposes
can be straightforwardly formulated in the context of simplicial
manifolds. Let $G$ be a simplicial Lie group. We denote the
corresponding face and degeneracy mappings $\partial_i$ and $s_i$,
respectively.
\begin{definition}\label{Moore}The \emph{Moore complex} $NG$
of $G$ is the Lie group chain complex $(NG,\delta)$ with
$$NG_n := \bigcap_{i=1}^{n} \ker \partial_i$$ and the differentials
$\delta_n: NG_n\to NG_{n-1}$ induced from the respective $0$th face
maps $\partial_0$ by restriction. It is a normal complex, i.e.
$\delta_n NG_n$ is a normal subgroup of  $NG_{n-1}$\footnote{It is a
normal subgroup of $G_{n-1}$ too.}. Of course, $NG_0=G_0$. The Moore
complex $NG$ is said to be of length $k$ if $NG_n$ is trivial for
$n>k$.\footnote{The objects of the full subcategory of simplicial
groups with Moore complex of length $k$ are also called
$k$\emph{-hypergroupoids} of groups \cite{Glenn}.}
\end{definition}

The Moore complex $NG$ carries a structure of a Lie
\emph{hypercrossed complex} structure, form which it can be
reconstructed \cite{Cond, Car}. To describe the idea behind this, we
will need following lemma.

\begin{lemma} Let $G$ be a simplicial Lie group. Then $G_n$ can be
decomposed as a semidirect product of Lie groups
$$G_n \cong \ker \partial_n \rtimes
s_{n-1}G_{n-1}.$$ Explicitly, for $g\in G_n$, the isomorphism is
given by
$$g\mapsto (g s_{n-1}\partial_n g^{-1}, s_{n-1}\partial_n g).$$
\end{lemma}

The following proposition \cite{Cond} is the result of a repetitive
application on the above lemma.
\begin{proposition}\label{reco}
For a simplicial Lie group $G$,
\begin{align*} G_n &\cong
(\ldots (NG_n\rtimes s_{0}NG_{n-1}))\rtimes \dots \rtimes
s_{n-1}\ldots s_0NG_1)
\end{align*}
\end{proposition}

The bracketing an ordering of the terms should be clear from the
first few terms of the sequence:
\begin{align}\label{decompAlt} G_1 &\cong
NG^0_1\rtimes s_{0}NG^0_{0}\nonumber\\
G_2 &\cong (NG^0_2\rtimes s_{0}NG^0_{1})\rtimes (s_1NG^0_1
\rtimes s_{1} s_0NG^0_0)\nonumber\\
G_3 &\cong ((NG^0_3\rtimes s_{0}NG^0_{2})\rtimes (s_1NG^0_2 \rtimes
s_{1} s_0NG^0_1))\rtimes\nonumber\\& ((s_{2}NG^0_{2}\rtimes
s_2s_0NG^0_1) \rtimes (s_{2} s_1NG^0_1 \rtimes s_2s_1s_0NG^0_0)).
\end{align}
We are not going to spell out the rather complicated definition of a
hypercrossed complex \cite{Car}. Instead, we give some examples.

\begin{example}\label{cm} A 1-hypercrossed complex of Lie groups is the same
thing as a Lie crossed module.
\end{example}

\begin{definition} Let $H$ and $D$ be two Lie
groups. We say that $H$ is a crossed $D$-module if there is a Lie
group morphism $\delta_1: H\to D$ and a smooth action of $D$ on $H$
$(d,h)\mapsto \,^d\hskip-0cm h$ such that
$$^{\delta_1 (h)}h' = hh'h^{-1}\,\, \mbox{(Peiffer condition)}$$ for $h,h'\in H$, and
$$\delta_1(^dh)= d \delta_1(h)d^{-1}$$ for $h\in H,
d\in D$ hold true.
\end{definition}

We will use the notation $H\stackrel{\delta_1}{\to} D$ or $H\to D$
for a crossed module.

\begin{definition}\label{CMM} A morphism between Lie crossed modules
$H\stackrel{\delta_1}{\to} D$ and $H'\stackrel{\delta_1'}{\to} D'$
is a pair of  Lie group morphisms $\lambda : H \to H'$ and $ \kappa
: D \to D'$ such that the diagram
$$
\begin{CD}
 H@>
 {\delta_1}>> D\\
 @V {\lambda}
 VV @V VV{\scriptstyle\kappa}\\
H'@>
 {\delta_1'}>>D'\,
\end{CD}
$$
commutes, and for any $h \in H$ and $d \in D$ we have the following
identity
\[ \lambda(^dh) = \,^{\kappa(d)}\lambda(h). \]
\end{definition}

Starting from a Lie crossed module $H\to D$ we can build up the
corresponding simplicial Lie group. Explicitly, cf. Proposition
\ref{reco},
$$G_0 = D, \hskip0.4cm G_1=(H\rtimes D), \hskip0.4cm  G_2 =
(H\rtimes(H\rtimes D)),\hskip0.4cm \mathrm{etc.}$$ The construction
can be interpreted as the internal nerve of the associated internal
category in the category of Lie groups (a strict Lie 2-group).
%%%%%%%%%%%%%%%%%%%%%%%%%%%%%%%%%%%%%%%%%%%%%%%%%%%%%%%%%%%%%%%%%%%%%%%
\begin{example} A Lie 2-hypercrossed complex is the
same thing as a Lie 2-crossed module \cite{Cond}.
\end{example}
\begin{definition} \label{2CM} A Lie 2-crossed module is a complex of Lie
groups
\begin{equation}\label{2crossmod}
H\stackrel{\delta_2}{\to} D\stackrel{\delta_1}{\to} K
\end{equation}
together with smooth left actions by automorphisms of $K$ on $H$ and
$D$ (and on $K$ by conjugation), and the Peiffer pairing, which is
an smooth equivariant map $\{\,,\,\}: D\times D \to H$, i.e.,
$^k\{d_1,
d_2\}=\{^k d_1, \,^k d_2\}$ such that:\\

i) (\ref{2crossmod}) is a complex of $K$-modules, i.e., $\delta_2$
and $\delta_1$ are $K$-equivariant and $\delta_2\delta_1 (h) =1$ for
$h\in H$,

ii) $d_1d_2d_1^{-1}=\delta_2\{d_1,d_2\}\,^{\delta_1(d_1)}d_2$, for
$d_1,d_2\in D$,

iii) $h_1h_2h_1^{-1}h_2^{-1}=\{\delta_2 h_1, \delta_2 h_2\}$, for
$h_1,h_2\in H$,

iv) $\{d_1 d_2, d_3\} = \{d_1, d_2 d_3 d_2^{-1}\}
\,^{\delta_1(d_1)}\{d_2, d_3\}$, for $d_1, d_2, d_3 \in D$,

v) $\{d_1, d_2 d_3\} = \,^{d_1 d_2 d_1^{-1}}\{d_1, d_3\}\{d_1,
d_2\}$, for $d_1, d_2, d_3 \in D$,

vi) $\{\delta_2(h), d\}\{d, \delta_2(h)\}= h^{\delta_1(d)}(h^{-1})$,
for
$d\in D, h\in H$,\\

\noindent wherein the notation $^kd$ and $^kh$ for left actions of
the element $k\in K$ on elements $d\in D$ and $h\in H$ has been
used.
\end{definition}

There is an obvious notion of a morphism of Lie 2-crossed modules.
\begin{definition}\label{2CMM} A morphism between Lie 2-crossed modules
$H\stackrel{\delta_2}{\to} D\stackrel{\delta_1}{\to} K$ and
$H'\stackrel{\delta_2'} {\to} D'\stackrel{\delta_1'}{\to}K'$ is a
triple of smooth group morphisms $H\to H'$, $D\to D'$ and $K\to K'$
making up, together with the maps $\delta_2$, $\delta_2'$,
$\delta_1$ and $\delta_1'$, a commutative diagram
\begin{equation}\label{mor2crossmod}
\begin{CD}
 H@>
 {\delta_2}>> D@>
 {\delta_1}>> K\\
 @V {\lambda}
 VV @V \mu VV @V VV{\scriptstyle\nu} \\
H'@>
 {\delta_2'}>>D'@>
 {\delta_1'}>>K'\,
\end{CD}
\end{equation}
and being compatible with the actions of $K$ on $D$ and $H$ and of
$K'$ on $D'$ and $H'$, respectively and with the respective Peiffer
pairings.
\end{definition}

The corresponding simplicial Lie group is given explicitly by, cf.
Proposition \ref{reco},
\begin{align*} &G_0 = K, \hskip0.4cm G_1=(D\rtimes K),
\hskip0.4cm  G_2 = (H\rtimes D)\rtimes(D\rtimes
K)),\\&G_3=(H\rtimes(H\rtimes D))\rtimes((H\rtimes
D)\rtimes(D\rtimes K)),\hskip0.4cm \mathrm{etc.}\end{align*}

This can be interpreted as an internal Duskin nerve \cite{Dus}.

\begin{example} A Lie 3-hypercrossed complex is the same thing as a
Lie 3-crossed module of \cite{Arvasi}.
\end{example}

We refer the interested reader to \cite{Car} for a thorough
discussion of hypercrossed complexes of groups and their relation to
simplicial groups.

At each level $n$, there is an lexicographically ordered set $S(n)$
of $2^n$ sets, which defines the compositions of the degeneracy maps
appearing in the decomposition of $G_n$. Explicitly for $S(n)$ we
have:
$$\{\emptyset< \{0\}<\{1\}< \{1,0\}< \{2\}<\{2,0\}< \{2,1\}<
\{2,1,0\}<\ldots< \{n-1,\ldots, 1\}< \{n-1,\ldots, 0\}\}.$$

The important role in the theory of hypercrossed complexes is played
by the actions $G_0\times NG_n \to NG_n$ defined by
$$g_0\times g_n\mapsto\,^{g_0}g_n:
(s_{n-1}\ldots s_0 g_0) g_n (s_{n-1}\ldots s_0 g_0)^{-1}$$ and the
so called \emph{Peiffer pairings}. In order to define these, we will
use the multi-indices like $\alpha$ and $\beta$ from $\bigcup_n
S(n)$ to write $s_\alpha$ for products of degeneracy maps

$$s_0,\,\, s_1,\,\,s_1s_0,\,\,s_2,\,\, s_2s_0,
\,\,s_2s_1,\,\,s_2s_1s_0,\ldots$$ In particular, for $g\in
NG_{n-\sharp\alpha}$ we have $s_\alpha g \in G_n$. For each $n$
consider the set $P(n)$ of pairs  $(\alpha, \beta)$ such that
$\emptyset<\alpha < \beta$ and $\alpha\cap\beta= \emptyset$, where
$\alpha\cap\beta$ is the set of indices belonging to both $\alpha$
and $\beta$.

\begin{definition}The Peiffer pairing (or lifting) $F_{\alpha,\beta}(g,h)\in NG_n$ for
$g\in NG_{n-\sharp\alpha}$, $h\in NG_{n-\sharp\beta}$ and $(\alpha,
\beta)\in P(n)$ is defined by
$$F_{\alpha,\beta}(g,h) = p_n(s_\alpha (g) s_\beta (h) s_\alpha (g)^{-1}s_\beta (h)^{-1}),$$
where $p_n$ is the projection to $NG_n$. For the projector $p_n$, we
have $p_n=p_n^1\ldots p^n_n$ with $p^i_n(g) = gs_{i-1}\partial_i
g^{-1}$.
\end{definition}

For us, the relevant Peiffer pairings at each level $n$ will be
those defined for pairs $(\alpha,\beta)\in P(n)$ such that
$\alpha\cup\beta = \{0,\ldots n\}$. We shall denote the set of such
pairs $\bar P(n)$.

\begin{remark} For a simplicial Lie algebra $\mathfrak{g}$, we have the
corresponding Moore complex $N\mathfrak{g}$ of Lie algebras, which
carries a structure of a hypercrossed complex of Lie algebras, cf.
\cite{Akca}. All the definitions and statements of this section
have, of course, their infinitesimal counterparts. Since these are
obvious, we shall not formulate them explicitly.

As shown by Quillen \cite{Quillen} there is an  adjunction between
simplicial Lie algebras and dg-Lie algebras. The part of the
adjunction going from simplicial Lie algebras to dg-Lie algebras
acts on the underlying simplicial vector spaces as the Moore complex
functor $N$.
\end{remark}

\section{Simplicial principal bundles}
Let $G$ be a simplicial Lie group and $X$ a simplicial manifold. In
this paper we use the name \emph{principal $G$-bundle} for a
\emph{twisted Cartesian product}. Therefore, we start with defining
twisting functions. Again, we will denote by $\partial_i$ and $s_i$
the corresponding face and degeneracy maps. We follow
\cite{May}.\footnote{Again, passing from sets to manifolds is
straightforward.}

\begin{definition}\label{tau} For a smooth function $\tau: X_n \to G_{n-1}$ to be a twisting,
the following conditions should be fulfilled:
$$ \partial_0\tau(x)\tau(\partial_0
x)= \tau(\partial_1x),$$
$$\partial_i\tau(x)=\tau(\partial_{i+1} x)\hskip0.4cm {\rm for}
\hskip0.4cm i>0,$$
$$s_i\tau(x)= \tau(s_{i+1} x)\hskip0.4cm {\rm for}
\hskip0.4cm i\geq0,$$
$$ \tau(s_0 x)=e_n\hskip0.4cm {\rm for}
\hskip0.4cm x \in X_n.$$
\end{definition}

\begin{definition} Let $\tau$ be a twisting function. A
twisted Cartesian product $P(\tau)= G\times_{\tau}X$ (alternatively
a principal $G$-bundle, or simply $G$-bundle, $P\to X$) is the
simplicial manifold with simplices
$$ P(\tau)_n= G_n\times X_n$$ and with the following face and degeneracy maps
$$ \partial_i(g,x)= (\partial_i g,\partial_i x) \hskip0.4cm {\rm for}
\hskip0.4cm i>0,$$
$$\partial_0(g,x)= (\partial_0 g
.\tau(x),\partial_0 x), $$
$$ s_i(g,x)= (s_i g, s_i x) \hskip0.4cm {\rm for} \hskip0.4cm i\geq
0.$$ The principal (left) $G$-action
$$G_n \times P(\tau)_n \to P(\tau)_n, \hskip0.4cm g'_n
\times (g_n, x_n)\mapsto (g'_n g_n, x_n)$$ and the projection
$$\pi_n:
P_n\to X_n, \hskip0.4cm (g_n, x_n)\mapsto x_n $$ are smooth
simplicial maps.
\end{definition}

Equivalence of two $G$-bundles $P(\tau)$ and $P(\tau')$ over the
same $X$  is described in terms of twisting as follows.
\begin{definition}\label{equivaltau} We call two twistings $\tau'$ and $\tau$
equivalent if there exists a smooth map $\psi: X \to G$ such that
$$
\partial_0\psi(x).\tau'(x)= \tau(x).\psi(\partial_0x),$$
$$ \partial_i\psi(x)=\psi(\partial_i x) \hskip0.4cm {\rm if}\hskip0.4cm i>0,$$
$$s_i\psi(x)=\psi(s_i x) \hskip0.4cm {\rm if}\hskip0.4cm i\geq 0.$$
In particular a twisting or the corresponding $G$-bundle $P(\tau)$
is trivial iff
$$\tau(x)= \partial_0\psi(x)^{-1}.\psi(\partial_0x),$$
with $\psi$ as above.
\end{definition}

As with ordinary bundles, simplicial principal bundles can be pulled
back and their structure groups can be changed using simplicial Lie
group morphisms. Twistings transform under these operations in an
obvious way.

There is a canonical construction of the classifying space
$\overline{W}G$ and of the universal $G$-bundle $WG$.

\begin{definition}\label{USimBun} The classifying space
$\overline{W}G$ is defined as follows. $\overline{W}G_0$ has one
element $\ast$ and $\overline{W}G_n=G_{n-1}\times G_{n-2}\times
\ldots \times G_0$ for $n>0$. Face and degeneracy maps are
$$s_0(\ast) = e_0, \hskip0.4cm  \hskip0.4cm \partial_i(g_0)=\ast
\hskip0.4cm {\rm for}\,\,\,\,i=0\,\,{\rm or}\,\,1$$
and
$$\partial_0 (g_{n}, \ldots g_0) =(g_{n-1},\ldots, g_0),$$
$$\partial_{i+1} (g_{n}, \ldots, g_0) =(\partial_ig_{n},\ldots ,
\partial_1g_{n-i+1},\partial_0g_{n-i}.g_{n-i-1}, g_{n-i-2} ,\ldots,  g_0),$$
$$s_0(g_{n-1}, \ldots ,g_0)=(e_n, g_{n-1},\ldots, g_0),$$
$$s_{i+1}(g_{n-1}, \ldots ,g_0)=(s_ig_{n-1},\ldots, s_0g_{n-i},e_{n-i},g_{n-i-1},
\ldots, g_0),$$ for $n>0$. With the choice of a twisting given by
$$\tau(g_{n-1}, \ldots ,g_0)=g_{n-1}$$ we have the universal $G$-principal
bundle
$$WG = G\times_\tau \overline{W}G.$$
\end{definition}

We have a relation between twistings and simplicial maps $X\to
\overline WG$ given by the following proposition.

\begin{proposition}\label{lemma}
The map $f_\tau: X \to \overline WG$ given by
$$x \mapsto \ast \hskip0.4cm  for\hskip0.4cm x\in X_0$$
and
$$x \mapsto (\tau(x), \tau (\partial_0 x),\ldots, \tau (\partial^i_0
x),\ldots, \tau (\partial^{n-1}_0 x))\hskip0.4cm  for\hskip0.4cm
x\in X_n,\, n>0$$ is a smooth simplicial map.

Vice versa, a smooth simplicial map $f: X \to \overline{W}_G$, given
by
$$x\mapsto \ast \hskip0.4cm for \hskip0.4cm x\in
X_0$$ and
$$x\mapsto (g^{(n)}_{n-1}(x), \ldots, g^{(n)}_0(x)) \hskip0.4cm for \hskip0.4cm x\in X_{n},\,\, n>0$$
defines a twisting by
$$\tau_f(x)= g^{(n)}_{n-1}(x) \hskip0.4cm for \hskip0.4cm x\in X_n, \,\,n>0.$$ We
have $\tau_{f_\tau} = \tau$ and $f_{\tau_f}= f$.
\end{proposition}

The role of the universal bundle is the following.
\begin{theorem} The principal $G$-bundle $G\times_\tau X$ corresponding to the
twisting $\tau$ is obtained from the universal bundle $WG$ as a
pullback under the simplicial map $f_\tau$.
\end{theorem}

\section{$L_\infty$-algebroids as 1-jets of simplicial sets}
This section is completely based on \cite{Sev}, to which we also
refer for the proofs. We keep, maybe with an occasional exception,
the notation and terminology used there. Let $\mathrm{SSM}$ denotes
the category with objects being surjective submersions between
supermanifolds and morphisms commutative squares. Any surjective
submersion $M\to N$ gives a simplicial supermanifold $X$, the nerve
of the the groupoid $M\times_N M \rightrightarrows N$. Further, let
$\mathrm{SSM}_1$ denote the full subcategory of $\mathrm{SSM}$ with
objects $\mathbb{R}^{0|1}\times N\stackrel{\mathrm{pr}_2}{\to} N$,
where $N$ is running through all supermanifolds. Let
$\mathrm{SM}_{[1]}$ be the category of supermanifolds with a right
action of the supersemigroup
$\mathbf{Hom}(\mathbb{R}^{0|1},\mathbb{R}^{0|1})$. Put in other
words, $\mathrm{SM}_{[1]}$ is the category of differential
non-negatively graded supermanifolds. We have the following lemma

\begin{lemma} The category $\widehat{\mathrm{SSM}}_1$ of presheaves  on $\mathrm{SSM}_1$ and the
category $\widehat{\mathrm{SM}}_{[1]}$ of presheaves on
$\mathrm{SM}_{[1]}$ are equivalent.
\end{lemma}
\begin{remark}
The above lemma follows from the useful observation
$${\rm{Hom}}(\mathbb{R}^{0|1}\times N\to N, \mathbb{R}^{0|1}\times X\to X)\simeq
{\rm{Hom}}(N, X)\times \mathbf{Hom}(\mathbb{R}^{0|1},
\mathbb{R}^{0|1})(N).\label{EqiuvSheaves}$$ Which just says that the
object $\mathbb{R}^{0|1}\times X \in \mathrm{SSM}_1$ corresponds to
the object $X\times \mathbf{Hom}(\mathbb{R}^{0|1},\mathbb{R}^{0|1})$
in $\mathrm{SM}_{[1]}$.
\end{remark}

\begin{definition} Let $F$ be a presheaf on $\mathrm{SSM}$. Its
restriction to $\mathrm{SSM}_1$ is an object in
$\widehat{\mathrm{SSM}}_1$. The corresponding object $F_1$ in
$\widehat{\mathrm{SM}}_{[1]}$ is called the \emph{1-jet} of $F$.
\end{definition}

\begin{remark} The representable 1-jets are of particular interest, since
they are represented by differential non-negatively graded
supermanifolds \footnote{Differential graded supermanifolds, i.e.
$Q$-supermanifolds provide a natural framework for the
Batalin-Vilkovisky formalism \cite{AKSZ}.}. Hence, they can provide
us with interesting examples of those. If the $\mathbb{Z}_2$ is
given by the parity of the $\mathbb{Z}$-degree, which will be always
the case in our examples, then we have a differential non-negatively
graded manifold. Let us recall, that a finite-dimensional,
positively graded differential manifold is the same thing as an
$L_\infty$-algebra. If it is only a non-negatively graded one then
it could be, for good reasons explained in \cite{Sev}, referred to
as an $L_\infty$-algebroid, cf. also \cite{Vor} for a formal
definition.
\end{remark}

Particular examples of presheaves on $\mathrm{SSM}$ come from
simplicial supermanifolds. If $K$ is a simplicial supermanifold and
$X$ the nerve of the groupoid defined by the surjective submersion
$M\to N$, the the corresponding sheaf $F^K \in
\widehat{\mathrm{SSM}}$ is defined by
$$F^K(M\to N) = \mathrm{Hom}(X,K),$$
i.e. it associates with the surjective submersion $M\to N$ the set
of all simplicial maps $X\to K$.

In \cite{Sev}, also the following sufficient condition for the 1-jet
$F^K_1$ of $F^K$ to be representable is given.

\begin{theorem} Let $K$ be a simplicial supermanifold fulfilling
the Kan conditions, which is moreover $m$-truncated for some $m\in
\mathbb{N}$. Then the 1-jet $F^K_1$ is representable.\footnote{We
refer to the proof as well for explanation of the Kan conditions and
$m$-truncatedness to the Appendix of \cite{Sev}. These notions were,
for simplicial manifolds, first introduced in \cite{Hen}. If $G$ is
a simplicial Lie group then $G$, $\overline W G$ and $WG$ fulfill
the Kan conditions\cite {DTU}.}
\end{theorem}

Another construction described in \cite{Sev} is the so called
1-approximation of a presheaf $F\in \widehat{\mathrm{SSM}}$. The
restriction of $F$ to ${\mathrm{SSM}}_1$ admits a right adjoint, the
induction.

\begin{definition}
The presheaf $\mathrm{app}_1F \in \widehat{\mathrm{SSM}}$ is defined
by successively applying the restriction and induction functors to
$F\in \widehat{\mathrm{SSM}}$.
\end{definition}

\begin{proposition}
If the jet 1-jet $F_1$ is represented by the differential
non-negatively graded supermanifold $X_F$ then the sheaf
$\mathrm{app}_1F\in \widehat{\mathrm{SSM}}$ is given by
\begin{align*} \mathrm{app}_1F (M\to N)&= \{\mathrm{morphisms\,\, of\,\, dg
\,\,manifolds}\,\, T[1](M\to N)\to X_F\}\\
&= \{\mathrm{morphisms\,\, of\,\, dg \,\,algebras}\,\,
C^\infty(X_F)\to \Omega(M\to N)\},
\end{align*}
where $T[1](M\to N)$ is the shifted fibrewise tangent bundle of $M$
and $\Omega(M\to N)= C^\infty(T[1](M\to N))$ are the fibrewise
differential forms on $M$.

If $X_F$ is positively graded then it can be identified with an
$L_\infty$-algebra $L_F$ and we have
$$\mathrm{app}_1F (M\to N)=\{\mathrm{Maurer-Cartan\,\,elements\,\,
of}\,\, L_F \otimes \Omega(M\to N)\}.$$
\end{proposition}

\begin{example}\label{pairGrd} Consider the presheaf in $\widehat{\mathrm{SSM}}$
represented by $Y\to X$.  Its 1-jet is 1-representable  by
$T[1](Y\to X)$ (equipped with the canonical differential), the
shifted fibrewise tangent bundle. This is just a fibrewise version
of the following well known fact $\mathrm{Hom} (R^{0|1}\times N\to
N, Y\to\ast)= \mathbf{Hom}(R^{0|1}, Y)(N)\cong \mathrm{Hom}(N,
T[1]Y)$, i.e. that ``maps from $\mathbb{R}^{0|1}$ to $M$ are the
same things as 1-forms on $M$''.
\end{example}

\begin{example}\label{GroupBundle}Let $G$ be a Lie group
with Lie algebra $\mathfrak{g}$ and $M\to N$ a surjective
submersion. A $G$-descent data on $M\to N$, i.e a descent of a
trivial $G$-bundle on $M$ to a $G$-bundle on $N$, is a map $g:
M\times_N M\to G$ satisfying $g(x,x)=e$ and $g(x,y)g(y,z)=g(x,z)$
for $(x,y,z)\in M\times_N M\times_N M$. The $G$-descent is the same
thing as a groupoid morphism from $M\times_N M\to G$. Let us
consider the presheaf
$$F(M\to N)=\{ G-\mathrm{descent\,\,data\,\,on\,\,}M\to N\},$$
which is in the above notation $F^{\mathcal{N}G}$, with
$\mathcal{N}G$ the nerve of $G$.

By definition, we have for the 1-jet $F_1^{\mathcal{N}G}$
$$F(\mathbb{R}^{0|1}\times N\to N)=\{
G-\mathrm{descent\,\,data\,\,on\,\,} \mathbb{R}^{0|1}\times N\to
N\}.$$ Such a $G$-descent data is a map $g:\mathbb{R}^{0|1}\times
\mathbb{R}^{0|1} \to G^N$ satisfying the above descent (1-cocycle)
condition and is equivalent to a map $\bar g : \mathbb{R}^{0|1} \to
G^N$, such that $\bar g(0)=e$. The relation between maps $g$ and
$\bar g$ is \footnote{Notice that $\bar g$ is just a trivialization
of the 1-cocycle $g$ over $\mathbb{R}^{0|1}$.}
$$g(\theta_0, \theta_1)= \bar g(\theta_0)^{-1} \bar g(\theta_1
),$$
$$\bar g(\theta)=g(0,\theta).$$

One way to see what the dg manifold representing the 1-jet looks
like is the following. Instead of imposing the condition $\bar
g(0)=e$, we can consider arbitrary functions $\bar g(\theta)$ modulo
left multiplications with the constant ones. So what we have is the
shifted tangent bundle $T[1]G$ equipped with the canonical
differential induced from the de Rham differential on $G$ modulo the
left $G$-action. This observation immediately leads to the dg
algebra of functions on
$\mathfrak{g}[1]\cong\mathbf{Hom}(\mathbb{R}^{0|1}, G)/G$ -- the
wedge algebra of left invariant forms on $G$ with the de Rham
differential, which is just the the Chevalley-Eilenberg complex of
$\mathfrak{g}$.

Equivalently,  with an obvious abuse of notation, which we will
commit also in the rest of the paper, we note that we can write
$$\bar g(\theta)= e - a\theta,$$ with $a\in \mathfrak{g}^N[1]$.
Hence, the 1-jet $F_1$ is
represented by the shifted Lie algebra
$\mathfrak{g}[1]$.
The differential is computed from $$g(\theta_0, \theta_1)= 1 +
a(\theta_0-\theta_1) + \frac{1}{2} [a,a]\theta_0\theta_1$$ by
computing
$$- (da)\theta_1=(\delta_\epsilon \bar
g(\theta_1)-e)=\frac{d}{d\epsilon}(g(\theta_0 + \epsilon,\theta_1 +
\epsilon)-e)\mid_{\epsilon=\theta_0=0} =
-\frac{1}{2}[a,a]\theta_1.$$
Hence, the differential is
$$da =\frac{1}{2}[a,a].$$

Finally, the functor $\mathrm{app}_1 F$ associates to a surjective
submersion $M\to N$ the set of flat fibrewise
connections.\footnote{Vice versa, if $G$ and all the fibres are
1-connected flat fibrewise  connections give us $G$-descents.}
\end{example}

\begin{example}\label{CrossedModuleBundle}Let $H\stackrel{\delta_1}\to D$
 be a crossed module of Lie groups
with the induced crossed module of Lie algebras
$$\mathfrak{h}\stackrel{\delta_1}\to \mathfrak{d}$$ and $M\to N$ a
surjective submersion. An $H\to D$-descent data on $M\to N$, is an
$(H\to D)$-valued 1-cocycle on on the groupoid $Y=M\times_N M$. Such
1-cocycles describe bundle gerbes, similarly as transition functions
describe principal bundles \cite{ACJ}. More explicitly, we have a
pair of maps $(h,d)$, $d: Y_1 \to D$ and $h: Y_2 \to H$, such that
$$d(y_1)d(y_2)=\delta_1(h(y_1,
y_2))d(y_1\circ y_2), \hskip0.4cm  \mathrm{for} \hskip0.4cm
(y_1,y_2)\in Y_2, $$
$$h(y_1, y_2)h(y_1\circ y_2, y_3) =\,^{d(y_1)}h(y_2, y_3)h(y_1,y_2\circ y_3)\hskip0.4cm  \mathrm{for} \hskip0.4cm
(y_1,y_2,y_3)\in Y_3,$$ and
$$d(e_x)=e \hskip0.4cm \mathrm{and} \hskip0.4cm h(e_{s(y)}, y)=h(y,
e_{t(y)})=e.$$ Let us consider the presheaf
$$F(M\to N)=\{ (H\to D)-\mathrm{descent\,\,data\,\,on\,\,}M\to N\},$$
which is in the above notation $F^{\mathcal{N}(H\to D)}$, with
$\mathcal{N}(H\to D)$ the Duskin nerve of $H\to D$.

By definition, we have for the 1-jet $F_1^{\mathcal{N}(H\to D)}$
$$F(\mathbb{R}^{0|1}\times N\to N)=\{
(H\to D)-\mathrm{descent\,\,data\,\,on\,\,} \mathbb{R}^{0|1}\times
N\to N\}.$$ Such an $(H\to D)$-descent data is a pair of maps
$(h,d)$, $d:\mathbb{R}^{0|1}\times \mathbb{R}^{0|1} \to D^N$ and
$h:\mathbb{R}^{0|1}\times \mathbb{R}^{0|1} \times
\mathbb{R}^{0|1}\to H^N$ satisfying the 1-cocycle condition, and is
equivalent to a pair of maps $(\bar h, \bar d)$, $\bar
d:\mathbb{R}^{0|1} \to D^N$ and $\bar h: \mathbb{R}^{0|1} \times
\mathbb{R}^{0|1}\to H^N$ such that
$$\bar d(0)=e, \hskip0.4cm \bar h(\theta,\theta)=\bar h(0,\theta)=e.$$ The relation between pairs of maps
$(d,h)$ and $(\bar d, \bar h)$ is\footnote{Again, the pair $(\bar h,
\bar d)$ is just a trivialization of the 1-cocycle $(h,d)$.}
$$d(\theta_0, \theta_1)= \bar d(\theta_0)^{-1} \delta_1(\bar h (\theta_0, \theta_1))
\bar d(\theta_1),$$
$$^{\bar d(\theta_0)}h(\theta_0, \theta_1, \theta_2)=
\bar h(\theta_0,\theta_1)\bar h (\theta_1, \theta_2)) \bar
h(\theta_0, \theta_2))^{-1}$$ and
$$\bar d (\theta)= d(0,\theta),$$
$$\bar h (\theta_0, \theta_1)= h(0,\theta_0,\theta_1).$$

Obviously, we can write
$$\bar d(\theta)= e - a\theta,$$ with $a\in
\mathfrak{d}[1]$ and
$$\bar h(\theta_0, \theta_1)= e + b\theta_0\theta_1,$$
with $b\in \mathfrak{h}[2]$. Hence, the 1-jet $F_1$ is
represented by the graded vector space $\mathfrak{d}[1]
\oplus \mathfrak{h}[2]$.

The differential is computed in a complete analogy with Example
\ref{GroupBundle} using expressions
$$d(\theta_0, \theta_1)=e + a(\theta_0 -
\theta_1) + (\frac{1}{2} [a,a] + \delta_1 b)\theta_0\theta_1,$$
$$h(\theta_0, \theta_1, \theta_2) =
e + b(\theta_0 \theta_1 + \theta_1 \theta_2-\theta_0 \theta_2)-\,^a
b\,\theta_0\theta_1\theta_2.$$

The resulting differential is:
$$da =\frac{1}{2}[a,a] + \delta_1 b,$$
$$db =\,^{a}b.$$

Since we have a positively graded dg manifold, we can describe it as
an $L_\infty$-algebra. It is actually a DGLA with generators only in
lowest two degrees, i.e a strict Lie 2-algebra. The nonzero
components are $L_0=\mathfrak{d}$ and $L_{-1}=\mathfrak{h}$. The
differential is $\delta_1: \mathfrak{h}\to \mathfrak{d}$. The
bracket on $\mathfrak{d}$ is given by its own Lie bracket, and the
bracket between $\mathfrak{d}$ and $\mathfrak{h}$ is given by the by
the action of $\mathfrak{d}$ on $\mathfrak{h}$. Let us note that Lie
2-algebras are one to one to crossed modules of Lie groups, cf.
\cite{BaezCrans}.
\end{example}

Finally, the functor $\mathrm{app}_1 F$ associates to a surjective
submersion $M\to N$ the set of  $(\mathfrak{h}\to
\mathfrak{d})$-valued flat fibrewise connections.
\section{$L_\infty$-algebra of $\overline W G$}
In this section we generalize the Examples \ref{GroupBundle} and
\ref{CrossedModuleBundle} to the case of a $G$-descent, where $G$ is
a simplicial Lie group  with Moore complex of length $k$. The
associated simplicial Lie algebra will be denoted by $\mathfrak{g}$.
Examples \ref{GroupBundle} and \ref{CrossedModuleBundle} correspond
to $k=1$ and $k=2$ respectively. Let $M\to N$ be a surjective
submersion. We define a $G$-descent data on $M\to N$ as a $G$-valued
twisting on the nerve of the groupoid $N\times_M N$. We recall, cf.
Definition \ref{tau}, that for $\tau: X_n \to G_{n-1}$ to be a
twisting, the following conditions should be fulfilled:
$$ \partial_0\tau(x)\tau(\partial_0
x)= \tau(\partial_1x),$$
$$\partial_i\tau(x)=\tau(\partial_{i+1} x)\hskip0.4cm {\rm for}
\hskip0.4cm i>0,$$
$$s_i\tau(x)= \tau(s_{i+1} x)\hskip0.4cm {\rm for}
\hskip0.4cm i\geq0,$$
$$ \tau(s_0 x)=e_n\hskip0.4cm {\rm for}
\hskip0.4cm x \in X_n.$$

Let us consider the presheaf
$$F(M\to N)=\{ G-\mathrm{descent\,\,data\,\,on\,\,}M\to N\},$$
which is in the notation of the previous section $F^{\overline WG}$,
i.e. the sheaf associating with the surjective submersion $N\to M$
the set of all simplicial maps from the nerve of the groupoid
$N\times_M N$ to the classifying space $\overline W G$.

By definition, we have for the 1-jet $F^{\overline W G}_1$
$$F(\mathbb{R}^{0|1}\times N\to N)=\{
G-\mathrm{descent\,\,data\,\,on\,\,} \mathbb{R}^{0|1}\times N\to
N\}.$$ Such an $G$-descent data is described by a twisting
\footnote{From now on we will omit the annoying $N$ and assume it
everywhere implicitly.} $\tau:(\mathbb{R}^{0|1})^{n} \to G_{n-1}^N$
and is equivalent to a function $\psi:(\mathbb{R}^{0|1})^{n} \to
G_{n}^N$ such that
$$ \partial_i\psi(\theta_0,\ldots\theta_n)=
\psi(\theta_0,\ldots,\hat{\theta}_i,\ldots\theta_n) \hskip0.4cm {\rm
if}\hskip0.4cm i>0,$$
$$s_i\psi(\theta_0,\ldots\theta_n)=\psi(\theta_0,\ldots,{\theta_i},{\theta_i}
\ldots\theta_n) \hskip0.4cm {\rm if}\hskip0.4cm i\geq 0.$$ We have
the following relation between $\tau$ and $\psi$ \footnote{As
before, $\psi$ is a trivialization of $\tau$, cf. Definition
\ref{equivaltau}.}
$$\tau(\theta_0, \ldots \theta_n)=
\partial_0\psi(\theta_0,\ldots,\theta_n)^{-1}\psi(\theta_1,\ldots,\theta_n),$$
$$\psi(\theta_0,\ldots,\theta_n)=\tau(0,\theta_0,\ldots,\theta_n).$$
From the definition of $\psi$ it follows that
$$\psi(0,\theta_1,\ldots,\theta_n)=\tau(0,0,\theta_1,\ldots,\theta_n)=
\tau(s_0(0,\theta_1,\ldots,\theta_n))=e_n.$$ Therefore, we write
$$\psi(\theta_0,\ldots,\theta_n)= 1 - a(\theta_1,\ldots \theta_n)\theta_0,$$
with $a(\theta_1,\ldots \theta_n)\in \oplus_{i=0}^n{n\choose
i}\mathfrak{g_n}[i+1].$ The function $a$ fulfils the following
identities
\begin{align}\partial_i a(\theta_1,\ldots\theta_n)&=
a(\theta_1,\ldots,\hat{\theta}_i,\ldots\theta_n) \hskip0.4cm {\rm
if}\hskip0.4cm i>0,\label{simplOna1}\\
\partial_0 a(0,\theta_1,\ldots \theta_n)&= a(\theta_1,\ldots
\theta_n)\label{simplOna2}\\
s_i a(\theta_1,\ldots\theta_n)&=
a(\theta_1,\ldots,\theta_i,\theta_i,\ldots,\theta_n) \hskip0.4cm
{\rm
if}\hskip0.4cm i>0,\label{simplOna3}\\
s_0a(\theta_1,\ldots,\theta_n)&=
a(0,\theta_1,\ldots,\theta_n)\label{simplOna4}.
\end{align} In the above list, the only possibly not completely obvious one is the
the $\partial_0$ equation (\ref{simplOna2}). However, this one
follows from the $s_0$ equation (\ref{simplOna4}) by an application
of $\partial_0$. From $(\ref{simplOna1})$ we immediately see that
$$a^{n}\in N\mathfrak{g}_n[n+1],$$
for the top component $a^n$ of $a(\theta_1,\dots \theta_n)=
a^n\theta_1\ldots\theta_n + \ldots$.

To proceed further, it will be more convenient change the Grassmann
coordinates by  $\bar\theta_0 = \theta_1$ and $\bar\theta_i=
\theta_{i+1} - \theta_{i}$ for $i>1$. In terms of $\bar\theta$s, we
get the following lemma for the decomposition of
$a(\bar\theta_0\ldots \bar\theta_{n-1})\in \oplus_{i=0}^{n}{n\choose
i}\mathfrak{g}_{n}[i+1]$ in terms of the shifted Moore complex
$N\mathfrak{g}_k[k+1]\oplus\ldots\oplus N\mathfrak{g}_0[1]$
\begin{lemma}\label{setminus} For $n\leq k$
$$a(\bar\theta_0,\ldots, \bar\theta_{n-1}) =
\sum_{\alpha\in S(n)}s_\alpha a^{n-\sharp
\alpha}\bar\theta^{S(n)\setminus \alpha},$$ where
$\bar\theta^\beta:=\bar\theta_{i_1}\ldots \bar\theta_{i_l}$ for
$\beta = \{{i_n},\ldots, i_1\} \in S(n)$.
\end{lemma}

\emph{Proof.} Straightforward computation using the fact that with
the new Grassmann variables $\bar \theta$ we have nice simplicial
relations $s_ia(\bar\theta_0,\ldots
\bar\theta_{n-1})=a(\bar\theta_0,\ldots, \bar\theta_{i-1}, 0,
\bar\theta_{i},\ldots, \bar\theta_{n-1})$ and $s_{i-1}\partial_i
a(\bar\theta_0,\ldots, \bar\theta_{n-1})= a(\bar\theta_0,\ldots,
\bar\theta_{n-1})|_{\bar\theta_{i-1}=0}\,\square.$

We see that, for $n\leq k$, the only independent component of $a \in
\oplus_{i=0}^n{n\choose i} \mathfrak{g_n}[i+1]$ is the top one
$a^{n}\in \mathfrak{g_n}[n+1]$.  Hence, the the 1-jet $F_1$ in this
case is represented by $ N\mathfrak{g}_k[k+1]\oplus\ldots\oplus
N\mathfrak{g}_0[1]$ as a graded manifold.

The differential can be obtained in analogy with the Examples
\ref{GroupBundle} and \ref{CrossedModuleBundle}. We write
\begin{align*}\tau(\theta_0,\ldots, \theta_n)&=
e+\partial_0a(\theta_1,\ldots,\theta_n)\theta_0-
a(\theta_2,\ldots,\theta_n)\theta_1 + \frac{1}{2}[\partial_0
a(\theta_1,\ldots
\theta_n),a(\theta_2,\ldots,\theta_n)]\theta_0\theta_1\\
&=e+\partial_0a(\theta_1,\ldots,\theta_n)\theta_0-
a(\theta_2,\ldots,\theta_n)\theta_1 + \frac{1}{2}[\partial_0
a(0,\theta_2\ldots
\theta_n),a(\theta_2,\ldots,\theta_n)]\theta_0\theta_1\\
&=e+\partial_0a(\theta_1,\ldots,\theta_n)\theta_0-
a(\theta_2,\ldots,\theta_n)\theta_1 + \frac{1}{2}[a(\theta_2,\ldots
\theta_n),a(\theta_2,\ldots,\theta_n)]\theta_0\theta_1.
\end{align*}

Now, using the above expression for $\tau$,
\begin{align*}&-(da(\theta_2,\ldots,
\theta_n))\theta_1=\delta_\epsilon(\psi(\theta_1,\ldots,
\theta_n)-e)= \frac{d}{d\epsilon}(\tau(\theta_0 +
\epsilon,\ldots,\theta_n + \epsilon)-e)\mid_{\epsilon=\theta_0=0}
\\&=\frac{d}{d\theta_1}\partial_0a(\theta_1,\ldots,
\theta_n)\theta_1+\sum_{i=2}^n \frac{d}{d\theta_i}a(\theta_2,\ldots,
\theta_n)\theta_1 -\frac{1}{2}[a(\theta_2,\ldots,
\theta_n),a(\theta_2,\ldots, \theta_n)]\theta_1.
\end{align*}
Hence, the differential is
\begin{align}\label{differential}da(\theta_1,\ldots,\theta_n) &=
-\frac{d}{d\theta_0}\partial_0a(\theta_0,\ldots,
\theta_n)-\sum_{i=1}^n \frac{d}{d\theta_i}a(\theta_1,\ldots,
\theta_n)\nonumber\\
&+\frac{1}{2}[a(\theta_1,\ldots, \theta_n),a(\theta_1,\ldots,
\theta_n)].
\end{align}

Now, we proceed in extracting the action $da^n$ of differential $d$
on the top component $a^n$. For this note: The first term gives
$-\partial_0 a^{n+1}$ and the second doesn't contribute to $da^n$ at
all . What remains is to determine the top component of the
commutator $[a(\theta_1,\ldots, \theta_n),a(\theta_1,\ldots,
\theta_n)]$. This leads to the following proposition.
\begin{proposition}\label{da} Let $G$ be a simplicial group with the simplicial
Lie algebra $\mathfrak{g}$. Assume that its Moore complex $NG$ is of
length $k$. Then the  1-jet $F_1$ of the simplicial manifold
$\overline W G$ is representable by the dg manifold $\oplus_{n=0}^k
N\mathfrak{g}_n[n+1]$. The differential $da^n$ on $a^n\in
N\mathfrak{g}_n[n+1]$ is described in terms of the face map
$\partial_0$, commutator of $\mathfrak{g}_0$, action of
$N\mathfrak{g}_0$ on $N\mathfrak{g}_n$ and Peiffer pairings
$f_{\alpha,\beta}$ with $(\alpha,\beta)\in \bar P(n)$ as follows:

For $n=0$
$$da^0 = -\partial_0 a^{1} + \frac{1}{2}[a^0,a^0],$$

for $n>0$
$$da^n = -\partial_0 a^{n+1} + \,^{a^0}
a^n + \sum_{(\alpha, \beta)\in\bar P(n)} \pm
f_{\alpha,\beta}(a^{n-\sharp\alpha},a^{n-\sharp\beta}),$$ where the
sign is given by the product of parity of
$(n-\sharp\alpha)(n-\sharp\beta +1)$ and the parity of the shuffle
defined by the pair $(S(n)\setminus\alpha,S(n)\setminus\beta)$.
\end{proposition}

\emph{Proof.} The 0th component is clear. What is left is to justify
the form of the second and third term in the above expression for
$da^n$, $n>0$. However, this is easily done using the above Lemma
\ref{setminus}. We just have to be careful about the degrees and
signs. We have
$$\frac{1}{2}[a(\bar\theta_0,\ldots,
\bar\theta_{n-1}),a(\bar\theta_0,\ldots, \bar\theta_{n-1})]^n=
[s_{n_1}\ldots s_0 a^0, a^n] + \sum_{(\alpha, \beta)\in\bar P(n)}\pm
[s_\alpha a^{n-\sharp\alpha},s_\beta a^{n-\sharp\beta}],$$ with the
sign given as the product of parities of
$(n-\sharp\alpha)(n-\sharp\beta +1)$ and of the shuffle defined by
the pair $(S(n)\setminus\alpha,S(n)\setminus\beta)$. Of course
$\sharp\alpha + \sharp\beta = n$.

The first term is just $^{a^0} a^n$, i.e. describing the action of
$N\mathfrak{g}_0$ on $N\mathfrak{g}_0$ shifted by 1 in degree.
Further, note that, by construction, the face $\partial_i$ and
degeneracy maps $s_i$ commute with the differential $d$. In
particular, it follows that $da^n$ must be in $\in
\bigcap_{i>0}\ker\partial_i[n+2]$. Moreover, since
$\partial_i\partial_0 =\partial_0\partial_{i+1}$, we also have
$\partial_0a^{n+1}\in \bigcap_{i>0}\ker\partial_i[n+2]$. Therefore,
we conclude that the sum over pairs $(\alpha,\beta)\in \bar P (n)$
in the above equation is also in $\in
\bigcap_{i>0}\ker\partial_i[n+2]$ and as such can be written,
trivially inserting the projection $p_n :\mathfrak{g}_n \to
N\mathfrak{g}_n$, as $\sum_{(\alpha, \beta)\in\bar P(n)}\pm
p_n[s_\alpha a^{n-\sharp\alpha},s_\beta a^{n-\sharp\beta}]$.
$\square$
\\

It is now straightforward to describe the $L_\infty$-algebra
corresponding to the above dg manifold explicitly. What we have is a
$k$-term DGLA $L=\oplus_{n=0}^{k}L_{-n}$ with components in degrees
$0,-1,\ldots -k$, given by $L_{-n}= N\mathfrak{g}_n$. The
differentials $d_n: N\mathfrak{g}_{n} \to N\mathfrak{g}_{n+1}$ are
given by the restrictions $d_n=\partial_0|_{N\mathfrak{g}_n}$ of the
zeroth face maps, i.e by the differentials $\delta_n$ of the Moore
complex $N\mathfrak{g}$, i.e, for $x_n\in{N\mathfrak{g}_n}$
\begin{equation}\label{delta1}
d_n x_n = \delta_n x_n .
\end{equation}
The only nonzero brackets are the binary brackets. The nonzero
binary brackets are determined by the following prescription:

The bracket $N\mathfrak{g}_0\times N\mathfrak{g}_0\to
N\mathfrak{g}_0$ is just the Lie bracket on $N\mathfrak{g}_0$, i.e
for $x_0\in N\mathfrak{g}_0$ and $y_0\in N\mathfrak{g}_0$
\begin{equation}\label{Lie}
[x_0, y_0].
\end{equation}

The brackets $N\mathfrak{g}_0\times N\mathfrak{g}_n\to
N\mathfrak{g}_n$: $(x, y)\mapsto [x_0,x_n] = - [x_n,x_0]$ are given
by the action of $N\mathfrak{g}_0$ on $N\mathfrak{g}_n$
\begin{equation}
[x_0,x_n]=- [x_n,x_0]=\,^{x_0}x_n.\label{action1}
\end{equation}

The bracket $N\mathfrak{g}_{n_1} \times N\mathfrak{g}_{n_2} \to
N\mathfrak{g}_{n}$ with $n=n_1+n_2$, for $n_1$ and $n_2$ nonzero, is
described as follows: For $x_{n_1}\in N\mathfrak{g}_{n_1}$ and
$x_{n_2}\in N\mathfrak{g}_{n_2}$
\begin{equation}\label{bracket1}[x_{n_1},x_{n_2}]=\sum_{(\alpha, \beta)\in\bar
P(n_1,n_2)} \pm f_{\alpha,\beta}(x_{n_1},x_{n_2}) +
(-1)^{(n_1+1)(n_2+1)}\sum_{(\alpha, \beta)\in\bar P(n_2,n_1)} \pm
f_{\alpha,\beta}(x_{n_2},x_{n_1})
\end{equation}
The $\pm$ sign is given by the product of parity of $n_1(n_2 +1)$
and the parity of the shuffle defined by the pair $(\alpha,\beta)\in
\bar P(n_1,n_2)$. Here $\bar P(n_1,n_2)\subset \bar P(n)$ denotes
the subset of $P(n)$ consisting of those pairs
$(\alpha,\beta)\in\bar P(n)$, for which $n-\sharp\alpha=n_1,
n-\sharp\beta=n_2$.

Let us now consider an arbitrary simplicial Lie algebra
$\mathfrak{g}$ with Moore complex of length $k$. Associated to
$\mathfrak{g}$ we have the (unique) simplicial group $G$ integrating
it, such that all its  components are simply connected. Therefore,
starting with $\mathfrak{g}$, we can consider the functor
$F^{\overline W G}_1$.  Correspondingly, we have the following
theorem:
\begin{proposition}\label{theMain}
Let $\mathfrak{g}$ be a simplicial Lie group with Moore complex
$N\mathfrak{g}$ of length $k$. Then $N\mathfrak{g}$ or becomes a
DGLA. The differential and the binary brackets are explicitly given
by formulas (\ref{delta1}-\ref{bracket1}). This DGLA structure on
$N\mathfrak{g}$ is the same one as described by Quillen's
construction in Proposition 4.4 of \cite{Quillen}.
\end{proposition}

We finish with some remarks:

\begin{remark}
Obviously, one can reformulate the above theorem in terms of a
$k$-hypercrossed complex of Lie algebras $\underline{\mathfrak{g}}$.
Such a $k$-hypercrossed complex $\underline{\mathfrak{g}}$ has a
structure of a $k$-term DGLA described by
(\ref{delta1}-\ref{bracket1}). \end{remark}

\begin{remark}
As noted above, crossed modules and Lie 2-algebras are one to one.
From the above theorem, we see that for $n>2$ only a part of the
full hypercrossed complex structure enters the description of the
DGLA $L$. For instance, already for $n=3$, only the symmetric part
of the Peiffer pairing appears. Nevertheless, the simplicial (Kan)
manifold $\overline W G$ can be interpreted as an integration of the
DGLA $N\mathfrak{g}$.

On the other hand, any $L_\infty$-algebra $L$, in particular any
DGLA, can be integrated to a (Kan) simplicial manifold $\int
 L$ \cite{Hen}. So one might try to compare the integration $\int N\mathfrak{g}$
with the Kan simplicial manifold $\overline W G$, or the
corresponding 1-jets (differentiations). Here we restrict ourselves
only to two related (obvious) remarks.

First, there is the following observation: Let $M$ be a simplicial
manifold. Assume that its corresponding 1-jet functor $F^M_1$ is
representable by an $L_\infty$ algebra $L$, with Chevalley-Eilenberg
complex $C(L)$. Also, let $\Omega({\Theta_N})$ be DGA of (normal)
forms on the simplicial supermanifold $\Theta_N$, the nerve
associated to the surjective submersion $\mathbb{R}^{0|1}\times N
\to N$.\footnote{See, e.g., \cite{Dupont} for the definition of
forms on simplicial sets.} Then, for the 1-jet corresponding to
$\int L$, we have $\Hom(\Theta_N, \int L)=\{\mathrm{morphisms\,\,
of\,\, dg \,\,algebras}\,\, C(L)\to \Omega(\Theta_N)\}$. This has to
be compared to 1-jet corresponding to $M$, i.e. to the set
$\Hom(\Theta_N, M)=\{\mathrm{morphisms\,\, of\,\, dg
\,\,algebras}\,\, C(L)\to \Omega(\mathbb{R}^{0|1}\times N \to N)\}$.

Second, in \cite{Hen}, simplicial homotopy groups
$\pi^{\textrm{spl}}_n \int L$ of $\int L$ have been shown to be
finite dimensional diffeological groups. Lie algebra of a
diffeological group is defined as the Lie algebra of its universal
cover, which is a Lie group. It is a result of \cite{Hen} that the
Lie algebra of $\pi^{\textrm{spl}}_n \int L$ is canonically
isomorphic to $H_{n-1}(L)$. Specified to the case of our interest:
the Lie algebra of $\pi^{\textrm{spl}}_n \int N\mathfrak{g}$ is
canonically isomorphic to  $H_{n-1}(N\mathfrak{g})$. If we now
consider $\overline WG$ just as a simplicial set then for its
simplicial homotopy groups we have $\pi_n \overline{W}G = \pi_{n-1}
G = H_{n-1}(NG)$. Given the Lie structure of $G$, $H_{n-1}(NG)$ and
hence  $\pi_n \overline{W}G$ can be considered as  diffeological
groups. In this sense we can talk about the Lie algebra of $\pi_n
\overline{W}G$, which is again $H_{n-1}(N\mathfrak{g})$.

\end{remark}

\begin{remark}\label{gauge1}
The functor $\mathrm{app}_1 F^{\overline{W}G}$ associates to a
surjective submersion $M\to N$ the set of $L$-valued flat  fibrewise
connections. To obtain also non-flat connections one may use the
Weil algebra of $L$, similarly as in \cite{SSS}. This, as well as
applications to higher gauge theory, will be described elsewhere.
\end{remark}

\begin{remark}\label{gauge2}
Regarding applications to higher gauge theory, in the forthcoming
work we plan to extend the results presented here to the case of
simplicial groupoids. In particular, we hope to describe a proper
generalization of the Atiyah groupoid to the simplicial case and
then obtain the Atiyah $L_\infty$-algebroid as a 1-jet of a properly
defined simplicial classifying space. This would lead us directly to
the notion of an $L$-valued connection on a simplicial principal
$G$-bundle (or on the corresponding hypercrossed complex bundle
gerbe). That this should be possible can be seen from the
description of connections and curvings of crossed module bundle
gerbes in \cite{Stev}, where a categorification of the Atiyah
algebroid is presented.
\end{remark}

\subsection*{Acknowledgments}
I would like to thank P. \v Severa for explanation of his work, D.
Stevenson for comments, T. Nikolaus for discussion and U. Schreiber
telling me about Quillen's construction. Also, I would like to thank
to IHES, LMU and MPIM for hospitality. The research was supported
under grant GA\v CR P201/12/G028.

\end{document}